\documentclass{amsart}
\usepackage{amssymb}

\usepackage{amscd}
\usepackage{thmdefs}

\input{tcilatex}
\theoremstyle{definition}
\theoremstyle{remark}
\numberwithin{equation}{section}

\input tcilatex

\begin{document}
\title[Moments of Block Operators]{Moments of the Block Operators in the Group Von Neumann Algebras}
\author{Ilwoo Cho}
\address[A. One and A. Two]{Author OneTwo address line 1\\
Author OneTwo address line 2}
\email[A. One]{aone@@aoneinst.edu}
\thanks{}
\date{}
\dedicatory{}
\thanks{}
\maketitle

\begin{abstract}
In this paper, we will consider the moments of the block operators of the
given group von Neumann algebra $L(G)$, where the given group $G$ is a
finitely presented discrete group $<X\,:\,R>$, where $X$ is the generator
set of $G$ and $R$ is the relation on the set $X,$ as the set of relators.
Define the canonical trace $tr$\ on $L(G)$ and the $W^{*}$-probability space 
$\left( L(G),tr\right) $ which is our free probabilitic object of this
paper. Define the block operators $T_{x}$ by $T_{x}$ $=$ $x$ $+$ $x^{-1}$ of 
$L(G),$ where $x\in X.$ In this paper, we will compute the moments and the
moment series of $T_{x},$ for $x$ $\in $ $X.$ By the computation, we can get
that if $x_{1}$ and $x_{2}$ are generators of presented groups $%
<X_{1}:\,R_{1}>$ and $<X_{2}:\,R_{2}>,$ respectively, and (i) if there 
\textbf{is} $n$ $\in $ $\Bbb{N}$ such that $x_{1}^{n}$ $\in $ $R_{1}$ and $%
x_{2}^{n}$ $\in $ $R_{2},$ or (ii) if there is \textbf{no} $n_{1},$ $n_{2}$ $%
\in $ $\Bbb{N}$ such that $x_{1}^{n_{1}}$ $\in $ $R_{1}$ and $x_{2}^{n_{2}}$ 
$\in $ $R_{2},$ then the block operators $x_{1}+x_{1}^{-1}$ in $L(G_{1})$
and $x_{2}+x_{2}^{-1}$ in $L(G_{2})$ are identically distributed.
\end{abstract}

\strut

The group von Neumann algebras are studied recently by various authors.
Group von Neumann algebras are interesting objects in Operator Algebra and
Free Probability. In this paper, we will consider the moments of certain
operators in group von Neumann algebras, where the group is presented by a
finite generator set and a finite relation. We will take a presented group $%
<X\,:\,R>,$ where $X$ is the generator set and $R$ is the relation on the
group, as the nonempty set of relators. For instance, the symmetric group $%
S_{3}$ can be presented by its generator set $X_{S_{3}}$ and its relation $%
R_{S_{3}},$ where

\strut \strut

\begin{center}
$X_{S_{3}}=\{a,b\}$
\end{center}

and

\begin{center}
$R_{s_{3}}=\{a^{2},b^{3},(ab)^{2}\}.$
\end{center}

\strut \strut

Let $H=\,<X\,:\,R>$ be a presented group and let $L(H)$ be the group von
Neumann algebra generated by $H.$ i.e.,

$\strut $

\begin{center}
$L(H)=\overline{\lambda (H)}^{w}=\overline{\Bbb{C}[H]}^{w},$
\end{center}

\strut

where $\lambda $ is the left (unitary) representation. If $a\in L(H)$ is an
operator, then it has the Fourier expansion,

\strut

\begin{center}
$a=\underset{h\in H}{\sum }a_{h}h,$ \ for \ $a_{h}\in \Bbb{C}$ and $h\in H.$
\end{center}

\strut

In fact, the group element $h$ in the previous expansion are understood as
unitary operators $\lambda _{h}$ on the Hilbert space $l^{2}(H).$ Recall that

\strut

\begin{center}
$x^{*}=\underset{h\in H}{\sum }\overline{a_{h}}\,h^{-1},$
\end{center}

\strut

where $h^{-1}=\lambda _{h^{-1}}=\lambda _{h}^{*}=h^{*},$ on $l^{2}(H),$ for
all $h\in H.$ Define the canonical trace $tr$ on $L(H)$ by

\strut

\begin{center}
$tr\left( x\right) =tr\left( \underset{h\in H}{\sum }a_{h}h\right)
=a_{e_{H}},$ for all $x\in L(H)$
\end{center}

\strut

where $e_{H}$ is the identity of $H.$ Then we can have the $W^{*}$%
-probability space $\left( L(H),tr\right) .$ The main purpose of this paper
is to compute the moments of the block operators of $L(H).$ In order to do
that we observed the free monoid $X^{*}$ $=$ $\left( X\text{ }\cup \text{ }%
X^{-1}\right) ^{\prime }$ of the group $H$ and the corresponding
combinatorial forms in $X^{*}$ of the group elements in $G.$ (If $Y$ is an
arbitrary set, then $Y^{\prime }$ is the set of all free words in $Y,$ which
is called the free monoid of $Y.$ The elements in $Y^{\prime }$ are called
the combinatorial forms of $<Y>,$ where $<Y>$ is the group generated by $Y.$)

\strut

Let $G=\,<X:\,R>$ be a presented group with its generator set $X$ $=$ $%
\{x_{1},$ ..., $x_{N}\}$ and its relation $R$ $=$ $\{r_{1},$ ..., $r_{M}\}$.
In Chapter 1, we will consider the free monoid $X^{*}$ of the group $G$
defined by

$\strut $

\begin{center}
$X^{*}\overset{def}{=}\cup _{n=0}^{\infty }\left\{
x_{j_{1}}^{p_{1}}...x_{j_{n}}^{p_{n}}: 
\begin{array}{l}
(j_{1},...,j_{n})\in \{1,...,N\}^{n} \\ 
(p_{1},...,p_{n})\in \{1,-1\}^{n}
\end{array}
\right\} ,$
\end{center}

\strut

which is the set of all free words of the generator set $X$ and $X^{-1}.$
When $n=0,$ the corresponding word is the empty word $\emptyset .$ There
exists a monoid homomorphism $\pi $ from $X^{*}$ onto the given group $G.$
Notice that, for any $g\in G,$ there is a subset $\pi ^{-1}(g)$ in $X^{*}.$
The elements $w_{g}$ in $\pi ^{-1}(g)$ are called the combinatorial forms of 
$g$ $\in $ $G.$ We will use the word problem on $X^{*},$ by computing the
moments of the block operators.

\strut

In Chapter 2, we will compute the moments of the block operators $%
T_{x}=x+x^{-1},$ for $x\in X,$ in the $W^{*}$-probability space $\left(
L(G),tr\right) .$ In particular, we have that ;

\strut \strut

(1) Suppose that there is no relator $r_{t}\in R$ such that $x^{k}=r_{t},$
for all $k\in \Bbb{N}.$ Then

\strut

\begin{center}
$tr\left( T_{x}^{m}\right) =\left[ 
\begin{array}{l}
\,m \\ 
\frac{m}{2}
\end{array}
\right] ,$ for all $m\in \Bbb{N},$
\end{center}

\strut

where $\left[ 
\begin{array}{l}
\,m \\ 
\frac{m}{2}
\end{array}
\right] \overset{def}{=}\left\{ 
\begin{array}{ll}
\left( 
\begin{array}{l}
m \\ 
\frac{m}{2}
\end{array}
\right) =\frac{m!}{\left( \frac{m}{2}!\right) \left( \frac{m}{2}!\right) } & 
\text{if }m\text{ is even} \\ 
\,\,\,\,\,\,\,\,\,\,\,\,\,\,\,\,\,\,\,\,\,\,\,\,\,\,\,\,0 & \text{if }m\text{
is odd.}
\end{array}
\right. $

\strut

(2) Suppose that there exist $r_{t}\in R$ and $n_{x}\in \Bbb{N}$ such that $%
r_{t}=x^{n_{x}}.$ Then

\strut

\begin{center}
$\ tr\left( T_{x}^{m}\right) =\left\{ 
\begin{array}{ll}
\left[ 
\begin{array}{l}
\,m \\ 
\frac{m}{2}
\end{array}
\right] & \text{if }m<n_{j} \\ 
&  \\ 
\left( 2^{k_{1}}-\left[ 
\begin{array}{l}
\,k_{1} \\ 
\frac{k_{1}}{2}
\end{array}
\right] \right) +\left[ 
\begin{array}{l}
\,m \\ 
\frac{m}{2}
\end{array}
\right] & \text{if }m=k_{1}n_{j}+k_{2},
\end{array}
\right. $
\end{center}

\strut

where $k_{1}\in \Bbb{N}$ and $k_{2}\in \Bbb{N}\cup \{0\}$ such that $0\leq
k_{2}<n_{x}.$

\strut $\strut $

Notice that if $G_{1}$ and $G_{2}$ are finitely presented groups (not
necessarily distinct) and $x_{1}$ and $x_{2}$ are generators of $G_{1}$ and $%
G_{2},$ respectively, and if there exists $n\in \Bbb{N}$ such that $%
x_{1}^{n} $ and $x_{2}^{n}$ are relators of $G_{1}$ and $G_{2},$
respectively, then the block operators $x_{1}+x_{1}^{-1}$ in $L(G_{1})$ and $%
x_{2}+x_{2}^{-1}$ in $L(G_{2})$ are identically distributed. Also, if $x_{1}$
and $x_{2}$ have no relators $x_{1}^{n}$ and $x_{2}^{k},$ for $n,$ $k$ $\in $
$\Bbb{N},$ then they are also identically distributed.

\strut \strut

\strut \strut

\section{Preliminaries}

\strut

\strut

In this paper, we will compute the moments of certain operators on a group
von Neumann algebra with its canonical faithful normal trace. We will
restrict our interests to finitely presented group von Neumann algebras. Let 
$G=\,<X\,:\,R>$ be a presented group, where $X$ is the finite generator set
of the group $G$ and $R$ is the relation on $G,$ as the set of all relators.
Denote the corresponding group von Neumann algebra by $L(G).$ Then each
operator $a\in L(G)$ has its Fourier expansion

\strut

(1.1) $\ \ \ \ \ \ \ \ \ \ \ a=\underset{g\in G}{\sum }a_{g}g,$ \ for \ $%
a_{g}\in \Bbb{C}$.

\strut

Note that we can regard $g$ in (1.3) as $\lambda _{g},$ for all $g\in G,$
where $\lambda $ is the left regular representation. Remark that $%
g^{*}=g^{-1},$ for all $g\in G$, in $L(G),$ and hence each operator $g\in
L(G)$ is unitary. For the group von Neumann algebra $L(G),$ we can define
the canonical trace $tr$ by

\strut

(1.2) $\ \ \ \ \ \ \ \ \ tr(a)\overset{def}{=}tr\left( \underset{g\in G}{%
\sum }a_{g}g\right) =a_{e_{G}},$

\strut

for all $a\in L(G),$ where $e_{G}$ is the identity of the group $G.$

\strut

\begin{definition}
Let $G=\,<X:\,R>$ be the presented group and $L(G),$ the corresponding group
von Neumann \ algebra. The algebraic pair $\left( L(G),tr\right) $ is called
the presented group $W^{*}$-probability space, where $tr$ is the canonical
trace given in (1.2). The operators in $\left( L(G),tr\right) $ are called
the random variables. Let $a\in L(G).$ Then the $n$-th moments of $a$ is
defined by

\strut 

$\ \ \ \ \ \ \ \ \ \ \ \ \ \ \ \ \ \ \ tr\left( a^{n}\right) ,$ \ \ \ for
all \ \ \ $n\in \Bbb{N}$.
\end{definition}

\strut

Now, we have our free probabilistic objects of this paper.

\strut

\strut \strut

\section{Moments of The Block Operators in Group von Neumann Algebras}

\strut

\strut

Throughout this chapter, let $G=\,<X:\,R>$ be the fixed presented graph with
its generator set $X$ and its relation $R$

\strut

\begin{center}
$X=\{x_{1},...,x_{N}\}$
\end{center}

and

\begin{center}
$R=\{r_{1},...,r_{M}\},$
\end{center}

\strut

where $M,\,N\,\in \Bbb{N}.$ In this chapter, we will compute the moments of
block operators $T_{x}=x+x^{-1}$, for $x,x^{-1}\in X,$ for $j$ $=$ $1$ $%
,..., $ $N.$

\strut

\begin{definition}
Let $G=\,<X:\,R>$ be a presented group with its generator set $X$ $=$ $%
\{x_{1},$ ..., $x_{N}\}$ and its relation $R$ $=$ $\{r_{1}$, $...,$ $r_{M}\},
$ where $M,$ $N$ $\in $ $\Bbb{N}$ and $r_{1},$ $...,$ $r_{M}$ are relators,
as elements in $X^{*}$. The operators $T_{j}$ are block operators in the
group von Neumann algebra $L(G),$ if

\strut 

$\ \ \ \ \ \ \ \ \ \ \ \ T_{j}=x_{j}+x_{j}^{-1},$ \ for all \ $j=1,...,N.$
\end{definition}

\strut \strut \strut

\strut

\strut

\subsection{Moments of Block Operators}

\strut

\strut

\strut

Let $Y$ be an arbitrary set. Then we can define a set $Y^{\prime },$
consisting of all free words in $Y.$ This set $Y^{\prime }$ is called the
free set of $Y.$ Let $G=\,<X:\,R>$ be the finitely presented group with

\strut

\begin{center}
$X=\{x_{1},...,x_{N}\}$ \ and \ $R=\{r_{1},...,r_{M}\}.$
\end{center}

\strut

Define the set $X^{*}$ be the free monoid $\left( X\text{ }\cup \text{ }%
X^{-1}\right) ^{\prime }$ of the set $X$ $\cup $ $X^{-1}.$ Notice that,
there exists the surjective (monoid) homomorphism $\pi $ $:$ $X^{*}$ $%
\rightarrow $ $G$ and, for any group element $g$ in $G,$ there exist words $%
w_{g}$ in $X^{*}$ satisfying that $\pi (w_{g})$ $=$ $g$ in $G.$ It is easy
to see that a corresponding word $w_{g}$ of $g$ is not uniquely determined.
We say that such words $w_{g}$ $\in $ $\pi ^{-1}(g)$ of $g$ $\in $ $G$ are
combinatorial forms of $g.$ Let $w_{g}$ $=$ $%
x_{j_{1}}^{p_{1}}...x_{j_{n}}^{p_{n}}$ $\in $ $X^{*}$ be a combinatorial
form of $g$ $\in G,$ where

$\strut $

\begin{center}
$(j_{1},$ ..., $j_{n})$ $\in $ $\{1,$ ..., $N\}^{n}$ and $(p_{1},$ ..., $%
p_{n})$ $\in $ $\{\pm 1\}^{n}$.
\end{center}

\strut

For convenience, we denote the word $x_{j_{n}}^{-p_{n}}...x_{j_{1}}^{-p_{1}}$
by $w_{g}^{-1}.$ Notice that $w_{g}^{-1}$ $\in $ $\pi ^{-1}(g^{-1})$ in $%
X^{*}.$

\strut

Now, fix the generators $x_{j}\in X$ and the block operator $%
T_{j}=x_{j}+x_{j}^{-1}.$ Consider the $n$-th moments $tr\left(
T_{j}^{m}\right) $ of $T_{j},$ for all $m\in \Bbb{N}.$ It is easy to see that

\strut

(2.1) \ \ $T_{j}^{m}=\underset{(p_{1},...,p_{m})\in \{1,-1\}^{m}}{\sum }%
\left( x_{j}^{p_{1}}\cdot \cdot \cdot x_{j}^{p_{m}}\right) =\underset{%
(p_{1},...,p_{m})\in \{1,-1\}^{m}}{\sum }x_{j}^{\sum_{k=1}^{m}p_{k}},$

\strut

for all $m\in \Bbb{N}.$ Notice that each word $x_{j}^{\sum_{k=1}^{m}p_{k}}$ $%
=$ $x_{j}^{p_{1}}$ $\cdot \cdot \cdot $ $x_{j}^{p_{m}}$ is regarded as an
element in the free monoid $X^{*}.$ i.e., without loss of generality, we can
consider the summands of $T_{j}^{m}$ as elements in the free monoid $X^{*}.$

\strut

\begin{lemma}
Let $T_{j}=x_{j}+x_{j}^{-1}$ be a block operator of the generator $x_{j}\in
X.$ If there is no $n_{j}\in \Bbb{N}$ such that $x_{j}^{n_{j}}$ $\in $ $R,$
then

\strut 

\ \ \ \ \ \ \ \ \ \ $tr\left( T_{j}^{m}\right) =\left\{ 
\begin{array}{ll}
\left( 
\begin{array}{l}
m \\ 
\frac{m}{2}
\end{array}
\right)  & \text{if }m\in 2\Bbb{N} \\ 
&  \\ 
\,\,\,\,\,\,\,\,0 & \text{if }m\in 2\Bbb{N}-1,
\end{array}
\right. $

\strut 

where $\left( 
\begin{array}{l}
n \\ 
k
\end{array}
\right) =\frac{n!}{k!\,(n-k)!},$ for $n$, $k\in \Bbb{N}.$
\end{lemma}

\strut

\begin{proof}
Assume that there is no relator $r_{t}\in R$ and $n_{j}\in \Bbb{N}$ such
that $r_{t}=x_{j}^{n_{j}}$ in $R$ $\subset $ $X^{*}.$ Then, by (2.1), we
have that

\strut

$\ \ \ tr\left( T_{j}^{m}\right) =tr\left( \underset{(p_{1},...,p_{m})\in
\{1,-1\}^{m}}{\sum }x_{j}^{\sum_{k=1}^{m}p_{k}}\right) $

\strut

$\ \ \ \ \ \ \ =tr\left( \underset{(p_{1},...,p_{m}),\,\sum_{j=1}^{m}p_{j}=0%
}{\sum }x_{j}^{\sum_{k=1}^{m}p_{k}}\right) $

\strut

$\ \ \ \ \ \ \ =tr\left( \underset{(p_{1},...,p_{m}),\,\sum_{k=1}^{m}p_{j}=0%
}{\sum }e_{G}\right) $

\strut

$\ \ \ \ \ \ \ =\left| \{(p_{1},...,p_{m})\in
\{1,-1\}^{n}:\sum_{k=1}^{n}p_{j}=0\}\right| $

\strut

$\ \ \ \ \ \ \ =\left( 
\begin{array}{l}
\,m \\ 
\frac{m}{2}
\end{array}
\right) .$

\strut \strut

The last equality holds because, to make $\sum_{k=1}^{m}p_{k}=0,$ the same
number of $+1$'s and $-1$'s should be appeared in the sequence $%
(p_{1},...,p_{m}).$ It is easy to see that if $m$ is odd, then the set

\strut

$\ \ \ \ \ \ \ \{(p_{1},...,p_{m})\in \{1,-1\}^{n}:\sum_{k=1}^{n}p_{j}=0\}$

\strut

is empty. Therefore, all odd moments of $T_{j}$ vanish.
\end{proof}

\strut \strut \strut \strut \strut

Now, we assume that there is a relator $r_{t}\in R$ such that $%
r_{t}=x_{j}^{n_{j}},$ for $n_{j}\in \Bbb{N},$ as a free word in $X^{*}.$
Clearly, the relator $r_{t}$ is a combinatorial form of $e_{G}$ (i.e., $\pi
(r_{t})=e_{G}$) in the group $G$ and the length $\left| r_{t}\right| $ of $%
r_{t}$ in $X^{*}$ is $n_{j}.$ Also, notice that if $r_{t}\in R,$ then the
words $wr_{t}w^{-1}$ and $wr_{t}^{-1}w^{-1}$ are also combinatorial forms of 
$e_{G},$ for all words $w$ and $w$ in $X^{*}$. By (2.1), we have that

\strut

\begin{center}
$T_{j}^{m}=\underset{(p_{1},...,p_{m})\in \{1,-1\}^{m}}{\sum }%
x_{j}^{\sum_{k=1}^{m}p_{k}}$ .
\end{center}

\strut

If $m=n_{j},$ then

\strut

\ $T_{j}^{n_{j}}=\left( x_{j}^{n_{j}}+x_{j}^{-n_{j}}\right) +\underset{%
(p_{1},...,p_{n_{j}})\in \{1,-1\}^{n_{j}},\,(p_{1},...,p_{n_{j}})\neq (\pm
1,...,\pm 1)}{\sum }x_{j}^{\sum_{k=1}^{n_{j}}p_{k}}$

\strut

(2.2)

\ \ $\ \ \ \ \ =2e_{G}+\underset{(p_{1},...,p_{n_{j}})\in
\{1,-1\}^{n_{j}},\,(p_{1},...,p_{n_{j}})\neq (\pm 1,...,\pm 1)}{\sum }%
x_{j}^{\sum_{k=1}^{n_{j}}p_{k}}.$

\strut

By the above formula (2.2), we can get the following lemma ;

\strut

\begin{lemma}
Let $x_{j}\in X$ and assume that there exists $t\in \{1,...,M\}$ such that $%
r_{t}$ $=$ $x_{j}^{n_{j}},$ for $n_{j}$ $\in $ $\Bbb{N}.$ Then

\strut 

$\ \ \ \ \ \ \ \ \ \ \ tr\left( T_{j}^{m}\right) =\left\{ 
\begin{array}{ll}
\left( 
\begin{array}{l}
\,m \\ 
\frac{m}{2}
\end{array}
\right)  & \text{for all even }m<n_{j} \\ 
&  \\ 
\,\,\,\,\,\,0 & \text{for all odd }m<n_{j}
\end{array}
\right. $

\strut and

$\ \ \ \ \ \ \ tr\left( T_{j}^{n_{j}}\right) =\left\{ 
\begin{array}{lll}
2 &  & \text{if }n_{j}\in 2\Bbb{N}-1 \\ 
&  &  \\ 
2+\left( 
\begin{array}{l}
\,n_{j} \\ 
\frac{n_{j}}{2}
\end{array}
\right)  &  & \text{if }n_{j}\in 2\Bbb{N}\text{.}
\end{array}
\right. $
\end{lemma}

\strut

\begin{proof}
The first formula is trivial, by the previous lemma.

\strut

Suppose that $n_{j}$ is an odd number in $\Bbb{N}.$ Then, by (2.2), we have
that

\strut

$\ \ \ T_{j}^{n_{j}}=2e_{G}+\underset{(p_{1},...,p_{n_{j}})\in
\{1,-1\}^{n_{j}},\,(p_{1},...,p_{n_{j}})\neq (\pm 1,...,\pm 1)}{\sum }%
x_{j}^{\sum_{k=1}^{n_{j}}p_{k}}.$

\strut

Since $n_{j}$ is an odd number, we cannot find the sequence $%
(p_{1},...,p_{n_{j}})$ in $\{\pm 1\}^{n_{j}}$ satisfying that $%
\sum_{k=1}^{n_{j}}p_{k}$ $=$ $0.$ So, we cannot find the $e_{G}$-terms in
the summand

\strut

$\ \ \ \ \ \ \ \ \ \underset{(p_{1},...,p_{n_{j}})\in
\{1,-1\}^{n_{j}},\,(p_{1},...,p_{n_{j}})\neq (\pm 1,...,\pm 1)}{\sum }%
x_{j}^{\sum_{k=1}^{n_{j}}p_{k}}$

\strut

of $T_{j}^{n_{j}}.$ Thus if $n_{j}$ is an odd number, then $tr\left(
T_{j}^{n_{j}}\right) =2.$

\strut

Now, assume that $n_{j}$ is an even number in $\Bbb{N}.$ Then, again by
(2.2), we have that

\strut

$\ tr\left( T_{j}^{n_{j}}\right) =2+\underset{(p_{1},...,p_{n_{j}})\in
\{1,-1\}^{n_{j}},\,(p_{1},...,p_{n_{j}})\neq (\pm 1,...,\pm 1)}{\sum }%
tr\left( x_{j}^{\sum_{k=1}^{n_{j}}p_{k}}\right) $

\strut

$\ \ \ \ \ \ \ =2+\left| \{(p_{1},...,p_{n_{j}})\in \{\pm
1\}^{n_{j}}:\sum_{k=1}^{n_{j}}p_{k}=0\}\right| $

\strut

$\ \ \ \ \ \ \ =2+\left( 
\begin{array}{l}
\,\,n_{j} \\ 
\frac{n_{j}}{2}
\end{array}
\right) $.

\strut
\end{proof}

\strut

Now, suppose that $m>n_{j}.$ There are two cases ;

\strut

\ \ \ (i) $\ m=k_{1}n_{j}+k_{2},$ where $1\leq k_{2}<n_{j}$ \ \ \ \ or

\ \ \ (ii) $m=kn_{j},$ for some $k\in \Bbb{N}.$

\strut

\begin{lemma}
Let $x_{j}\in X$ and assume that there exists $t\in \{1,...,M\}$ such that $%
r_{t}$ $=$ $x_{j}^{n_{j}},$ for $n_{j}$ $\in $ $\Bbb{N}$. Then

\strut 

(1) If $m=kn_{j},$ for $k\in \Bbb{N},$ then

\strut 

$\ \ \ \ \ \ \ \ \ tr\left( T_{j}^{m}\right) =\left( 2^{k}-\left[ 
\begin{array}{l}
\,\,k \\ 
\frac{k}{2}
\end{array}
\right] \right) +\left[ 
\begin{array}{l}
\,kn_{j} \\ 
\frac{kn_{j}}{2}
\end{array}
\right] .$

\strut \strut \strut 

(2) If $m=k_{1}n_{j}+k_{2},$ for $k_{1},k_{2}\in \Bbb{N}$ and if $1\leq
k_{2}<n_{j},$ then

\strut 

$\ \ \ \ \ \ \ \ \ tr\left( T_{j}^{m}\right) =2^{k_{1}}+\left[ 
\begin{array}{l}
m \\ 
\frac{m}{2}
\end{array}
\right] -\left[ 
\begin{array}{l}
k_{1} \\ 
\frac{k_{1}}{2}
\end{array}
\right] ,$

\strut 

where

$\ \ \ \ \ \ \ \ \ \left[ 
\begin{array}{l}
\,\,t \\ 
\frac{t}{2}
\end{array}
\right] \overset{def}{=}\left\{ 
\begin{array}{lll}
\left( 
\begin{array}{l}
\,\,t \\ 
\frac{t}{2}
\end{array}
\right)  &  & \text{if }t\text{ is even} \\ 
&  &  \\ 
0 &  & \text{if }t\text{ is odd.}
\end{array}
\right. $
\end{lemma}

\strut

\begin{proof}
(\textbf{1}) Let $(p_{1},...,p_{m})\in \{\pm 1\}^{m},$ where $m=kn_{j}$ is
sufficiently big number in $\Bbb{N},$ where $k,n_{j}\in \Bbb{N}.$ Define
subsequences

\strut

$\ \ \ \ \ \ \ \mathbf{i}_{+}=\left( \underset{n_{j}\text{-times}}{%
\underbrace{1,.....,1}}\right) $ \ \ \ and \ \ \ $\mathbf{i}_{-}=\left( 
\underset{n_{j}\text{-times}}{\underbrace{-1,.......,-1}}\right) .$

\strut

Since $m=kn_{j},$ there exists a sequence $P=(p_{1},...,p_{m})$ such that

\strut

$\ \ \ \ \ \ \ \ \ P=\left( \mathbf{i}_{i_{1}},....,\mathbf{i}%
_{i_{k}}\right) ,$ for $i_{1},...,i_{k}\in \{+,\,-\}.$

\strut

We define the set $W_{j}$, consisting of such sequences. i.e.,

\strut

(2.3)$\ \ \ \ \ \ W_{j}\overset{def}{=}\{\left( \mathbf{i}_{i_{1}},...,%
\mathbf{i}_{i_{k}}\right) :i_{1},...,i_{k}\in \{+,\,-\}\}.$

\strut

Note that, for $\mathbf{i}_{+},$ we have $x_{j}^{n_{j}}=r_{t_{j}}$ and, for $%
\mathbf{i}_{-},$ we have $x_{j}^{-n_{j}}=r_{t_{j}}^{-1}$ in $X^{*}.$ Thus,
we can get that

\strut

(2.4) $\ \ \ \ \ \ \ \ \ \ \ \ \ \ \ \ \left| W_{j}\right| =2^{k}$ \ \ \ \ \
\ \ \ and

\strut

(2.5) $\ \ \ \ \ \ \ \ \sum_{k=1}^{m}p_{k}=pn_{j},$ where $p=1$ or ... or $%
k, $

\strut

for all $(p_{1},...,p_{m})$ in $W_{j}.$ We will define a subset $%
W_{j}^{\prime }$ of the set $W_{j}$ (if exists) by

\strut

\strut (2.6)

$\ W_{j}^{\prime }=\left\{ \left( \mathbf{i}_{i_{1}},...,\mathbf{i}%
_{i_{k}}\right) : 
\begin{array}{l}
i_{j_{1}}=...=i_{j_{\frac{k}{2}}}=+, \\ 
i_{s_{1}}=...=i_{s_{\frac{k}{2}}}=-, \\ 
\{i_{j_{1}},...,i_{j_{\frac{k}{2}}}\}\cup \{i_{s_{1}},...,i_{s_{\frac{k}{2}%
}}\}=\{i_{1},...,i_{k}\}
\end{array}
\right\} $

\strut

Then $W_{j}^{\prime }\subseteq W_{j}$ and, since there are same numbers of $%
+ $'s and $-$'s, we can have

\strut

(2.7) $\ \ \ \ \ \ \ \ \ \ \ \left| W_{j}^{\prime }\right| =\left( 
\begin{array}{l}
\,k \\ 
\frac{k}{2}
\end{array}
\right) ,$

\strut

if the nonempty subset $W_{j}^{\prime }$ exists in $W_{j}.$ (It is easily
see that if $k$ is even, then $W_{j}^{\prime }$ exists in $W_{j}.$ And if $k$
is odd, then $W_{j}^{\prime }$ is empty.) Define the subset $S_{0}$ of
sequences in $\{\pm 1\}^{m}$ by

\strut

\ $\ \ \ \ \ \ \ S_{0}=\{(p_{1},...,p_{m})\in \{\pm
1\}^{m}:\sum_{k=1}^{m}p_{k}=0\}.$

\strut

Then, by (2.6), $W_{j}^{\prime }\subset S_{0}.$ In fact,

\strut

(2.8) \ \ \ \ $S_{0}\cup W_{j}^{\prime }=S_{0}$ and $\left(
W_{j}\,\,\setminus \,\,W_{j}^{\prime }\right) \cap S_{0}=\emptyset .$

\strut

Assume that $k$ is even. Then

\strut

$\ \ tr\left( T_{j}^{m}\right) =tr\left( \underset{(p_{1},...,p_{m})\in
\{\pm 1\}^{m}}{\sum }x_{j}^{\sum_{k=1}^{m}p_{k}}\right) $

\strut

$\ \ \ \ \ \ \ =tr\left( \underset{(p_{1},...,p_{m})\in W_{j}\,\,\setminus
\,\,W_{j}^{\prime }}{\sum }x_{j}^{\sum_{k=1}^{m}p_{k}}\right) $

\strut

$\ \ \ \ \ \ \ \ \ \ \ \ \ \ \ \ +tr\left( \underset{(p_{1},...,p_{m})\in
\{\pm 1\}^{m},\,\,(p_{1},...,p_{m})\in S_{0}}{\sum }x_{j}^{%
\sum_{k=1}^{m}p_{k}}\right) $

\strut

$\ \ \ \ \ \ \ =\left| W_{j}\,\,\setminus \,\,W_{j}^{\prime }\right|
+tr\left( \underset{(p_{1},...,p_{m})\in \{\pm
1\}^{m},\,\,(p_{1},...,p_{m})\in S_{0}}{\sum }e_{G}\right) $

\strut

where $W_{j}$ is defined in (2.3), by (2.8)

\strut

$\ \ \ \ =\left( 2^{k}-\left( 
\begin{array}{l}
\,\,k \\ 
\frac{k}{2}
\end{array}
\right) \right) +\left| S_{0}\right| $

\strut

\strut (2.9)

$\ \ \ =\left\{ 
\begin{array}{ll}
\left( 2^{k}-\left( 
\begin{array}{l}
\,\,k \\ 
\frac{k}{2}
\end{array}
\right) \right) +\left( 
\begin{array}{l}
\,m \\ 
\frac{m}{2}
\end{array}
\right) & \text{if }m\text{ is even} \\ 
&  \\ 
2^{k}-\left( 
\begin{array}{l}
\,\,k \\ 
\frac{k}{2}
\end{array}
\right) & \text{if }m\text{ is odd}
\end{array}
\right. ,$

\strut

where $k\in 2\Bbb{N}$ and $m\in \Bbb{N}.$ Now, let's suppose that $k$ is an
odd number greater than 1 in $\Bbb{N}.$ Then we can have that

\strut

(2.10) $\ \ \ \ \ \ W_{j}^{\prime }=\emptyset $ \ \ and \ \ $W_{j}\cap
S_{0}=\emptyset .$

\strut

Therefore, by (2.10), we have that

\strut

$tr\left( T_{j}^{m}\right) =tr\left( \underset{(p_{1},...,p_{m})\in W_{j}}{%
\sum }x_{j}^{\sum_{k=1}^{m}p_{k}}\right) $

\strut

\ $\ \ \ \ \ \ \ \ \ \ \ \ \ \ +tr\left( \underset{(p_{1},...,p_{m})\in
\{\pm 1\}^{m},\,\,(p_{1},...,p_{m})\in S_{0}}{\sum }x_{j}^{%
\sum_{k=1}^{m}p_{k}}\right) $

\strut

\ \ $\ \ \ \ \ \ \ =\left| W_{j}\right| +tr\left( \underset{%
(p_{1},...,p_{m})\in \{\pm 1\}^{m},\,\,(p_{1},...,p_{m})\in S_{0}}{\sum }%
e_{G}\right) $

\strut

$\ \ \ \ \ \ \ \ \ =\left| W_{j}\right| +\left| S_{0}\right| $

\strut

\strut (2.11)

$\ \ \ \ \ \ \ \ \ \ =\left\{ 
\begin{array}{lll}
2^{k}+\left( 
\begin{array}{l}
\,m \\ 
\frac{m}{2}
\end{array}
\right) &  & \text{if }m\text{ is even} \\ 
&  &  \\ 
2^{k} &  & \text{if }m\text{ is odd.}
\end{array}
\right. $

\strut

Now, define a new notation

\strut

$\ \ \ \ \ \ \ \left[ 
\begin{array}{l}
\,\,t \\ 
\frac{t}{2}
\end{array}
\right] \overset{def}{=}\left\{ 
\begin{array}{lll}
\left( 
\begin{array}{l}
\,\,t \\ 
\frac{t}{2}
\end{array}
\right) &  & \text{if }t\in 2\Bbb{N} \\ 
&  &  \\ 
0 &  & \text{if }t\in 2\Bbb{N}-1.
\end{array}
\right. $

\strut

Then the formuli (2.9) and (2.11) can be shortened by

\strut

$\ \ \ \ \ \ \ \ \ tr\left( T_{j}^{kn_{j}}\right) =\left( 2^{k}-\left[ 
\begin{array}{l}
\,\,k \\ 
\frac{k}{2}
\end{array}
\right] \right) +\left[ 
\begin{array}{l}
\,kn_{j} \\ 
\frac{kn_{j}}{2}
\end{array}
\right] .$

\strut \strut

(\textbf{2}) Let's assume that $m=k_{1}n_{j}+k_{2},$ where $n_{j}\nmid k_{2}$
and $1\leq k_{2}<n_{j}$. By (2.1), we have that

\strut

$\ \ \ \ \ \ \ \ \ \ \ \ \ T_{j}^{m}=\underset{(p_{1},...,p_{m})\in \{\pm
1\}^{m}}{\sum }x^{\sum_{k=1}^{m}p_{k}}.$

\strut \strut

Let's regard the summands $x^{\sum_{k=1}^{m}p_{k}}$ as free words in the
free monoid $X^{*}.$ Then there is a set

\strut

$\ \ \ \ \ \ \ S_{0}=\{(p_{1},...,p_{m})\in \{\pm
1\}^{m}:\sum_{k=1}^{m}p_{k}=0\}$

\strut

with its cardinality $\left| S_{0}\right| =\left[ 
\begin{array}{l}
m \\ 
\frac{m}{2}
\end{array}
\right] .$ Now, define a set of free words $W^{j}$ by

\strut

$\ W^{j}=\left\{ \left( (p_{t_{1}},...,p_{t_{k_{2}}})\leadsto (\mathbf{i}%
_{i_{1}},...,\mathbf{i}_{i_{k_{1}}})\right) : 
\begin{array}{l}
p_{t_{k}}\in \{\pm 1\}, \\ 
\,i_{k}\in \{+,\,-\},
\end{array}
\right\} ,$

\strut

where $\mathbf{i}_{+}$ and $\mathbf{i}_{-}$ are defined in (1) and $\leadsto 
$ means the insertion. i.e., the sequence

\strut

$\ \ \ \ \ \ \ \ \ \ \ \ \ \ \ \left( (p_{t_{1}},...,p_{t_{k_{2}}})\leadsto (%
\mathbf{i}_{i_{1}},...,\mathbf{i}_{i_{k_{1}}})\right) $

\strut

is the free word in $X^{*}$ with its length $m$ $=$ $k_{1}n_{j}+k_{2}.$ For
example,

\strut

\ $\ \ \ \ \ \ \ \ \ \ \ \ \ \ \ \ \ \left( (p_{1})\leadsto (\mathbf{i}_{+},%
\mathbf{i}_{+})\right) $

is

$\ \ \ \ \left( \mathbf{i}_{+},p_{1},\mathbf{i}_{+}\right) $ \ or \ $\left(
p_{1},\mathbf{i}_{+},\mathbf{i}_{+}\right) $ \ or \ $\left( \mathbf{i}_{+},%
\mathbf{i}_{+},p_{1}\right) .$

\strut

where $m=2n_{j}+1.$ Define the subset $W_{k_{2}}^{j}$ of $W^{j}$ (if exists)
by

\strut

$\ W_{k_{2}}^{j}=\{\left( (p_{t_{1}},...,p_{t_{k_{2}}})\leadsto (\mathbf{i}%
_{i_{1}},...,\mathbf{i}_{i_{k_{1}}})\right)
:\sum_{i=1}^{k_{2}}p_{t_{i}}=0\}. $

\strut

Let's assume that $W_{k_{2}}^{j}$ exists. Then we can define the subset $%
W_{k_{2}}^{j}(S_{0})$ of $W_{k_{2}}^{j}$ by

\strut

$W_{k_{2}}^{j}(S_{0})=\left\{ \left( (p_{t_{1}},...,p_{t_{k_{2}}})\leadsto (%
\mathbf{i}_{i_{1}},...,\mathbf{i}_{i_{k_{1}}})\right) : 
\begin{array}{l}
\sum_{i=1}^{k_{2}}p_{t_{i}}=0 \\ 
(\mathbf{i}_{i_{1}},...,\mathbf{i}_{i_{k_{1}}})\in W_{j}^{\prime }
\end{array}
\right\} ,$

\strut

where $W_{j}^{\prime }$ is defined in (1). It is easy to see that if both $%
k_{1}$ and $k_{2}$ are even, then $W_{k_{2}}^{j}$ exists in $W_{k_{2}}^{j}.$
Then

\strut

(2.12) $\ \ \ \ \ \ \ \ \left| W_{k_{2}}^{j}\right| =2^{k_{1}}+\left[ 
\begin{array}{l}
k_{2} \\ 
\frac{k_{2}}{2}
\end{array}
\right] ,$

\strut

(2.13) \ \ \ \ \ \ \ $\left| W_{k_{2}}^{j}(S_{0})\right| =\left[ 
\begin{array}{l}
k_{1} \\ 
\frac{k_{1}}{2}
\end{array}
\right] +\left[ 
\begin{array}{l}
k_{2} \\ 
\frac{k_{2}}{2}
\end{array}
\right] $,

\strut

(2.14) \ \ \ \ \ \ \ \ \ $W_{k_{2}}^{j}(S_{0})\cup S_{0}=S_{0}$,

\strut

(2.15)\ \ $\left| S_{0}\,\,\,\setminus \,\,W_{k_{2}}^{j}(S_{0})\right|
=\left[ 
\begin{array}{l}
m \\ 
\frac{m}{2}
\end{array}
\right] -\left( \left[ 
\begin{array}{l}
k_{1} \\ 
\frac{k_{1}}{2}
\end{array}
\right] +\left[ 
\begin{array}{l}
k_{2} \\ 
\frac{k_{2}}{2}
\end{array}
\right] \right) $

\strut

By (2.12) and (2.15), we can compute that ;

\strut

$\ \ tr\left( T_{j}^{m}\right) =tr\left( \underset{(p_{1},...,p_{m})\in
\{\pm 1\}^{m}}{\sum }x^{\sum_{k=1}^{m}p_{m}}\right) $

\strut

$\ \ \ \ =tr\left( \underset{(p_{1},...,p_{n})\in W_{k_{2}}^{j}}{\sum }%
x^{\sum_{k=1}^{m}p_{m}}\right) $

\strut

$\ \ \ \ \ \ \ \ \ \ \ \ \ \ +tr\left( \underset{(p_{1},...,p_{m})\in
S_{0}\,\setminus \,\,W_{k_{2}}^{j}(S_{0})}{\sum }x^{\sum_{k=1}^{m}p_{m}}%
\right) $

\strut

$\ \ \ \ =\left| W_{k_{2}}^{j}\right| +\left| S_{0}\,\,\setminus
\,\,W_{k_{2}}^{j}(S_{0})\right| $

\strut

$\ \ \ \ =\left( 2^{k_{1}}+\left[ 
\begin{array}{l}
k_{2} \\ 
\frac{k_{2}}{2}
\end{array}
\right] \right) +\left( \left[ 
\begin{array}{l}
m \\ 
\frac{m}{2}
\end{array}
\right] -\left( \left[ 
\begin{array}{l}
k_{1} \\ 
\frac{k_{1}}{2}
\end{array}
\right] +\left[ 
\begin{array}{l}
k_{2} \\ 
\frac{k_{2}}{2}
\end{array}
\right] \right) \right) $

\strut \strut

$\ \ \ \ =2^{k_{1}}+\left[ 
\begin{array}{l}
m \\ 
\frac{m}{2}
\end{array}
\right] -\left[ 
\begin{array}{l}
k_{1} \\ 
\frac{k_{1}}{2}
\end{array}
\right] .$

\strut
\end{proof}

\strut \strut

Remark that, by the previous lemma, we have that if $k=1,$ then,

\strut

\begin{center}
$tr\left( T_{j}^{n_{j}}\right) =\left( 2^{1}-\left[ 
\begin{array}{l}
1 \\ 
\frac{1}{2}
\end{array}
\right] \right) +\left[ 
\begin{array}{l}
n_{j} \\ 
\frac{n_{j}}{2}
\end{array}
\right] =2+\left[ 
\begin{array}{l}
n_{j} \\ 
\frac{n_{j}}{2}
\end{array}
\right] .$
\end{center}

\strut \strut

By the previous lemmas, we can get the following theorem ;

\strut

\begin{theorem}
Let $G=\,<X:\,R>$ be a finitely presented group with its generator set $X$ $=
$ $\{x_{1},$ ..., $x_{N}\}$ and its relation $R$ $=$ $\{r_{1},$ ..., $%
r_{M}\}.$ Fix a generator $x_{j}\in X$ satisfying that $r_{t}=x_{j}^{n_{j}},$
for some $n_{j}\in \Bbb{N}\,\setminus \,\{1\},$ where $r_{t}\in R.$ Then

\strut 

(1) if $m<n_{j},$ then $tr\left( T_{j}^{m}\right) =\left[ 
\begin{array}{l}
m \\ 
\frac{m}{2}
\end{array}
\right] ,$

\strut 

(2) if $m=kn_{j},$ for $k\in \Bbb{N},$ then

\strut 

$\ \ \ \ \ \ \ \ \ \ \ tr\left( T_{j}^{m}\right) =\left( 2^{k}-\left[ 
\begin{array}{l}
\,\,k \\ 
\frac{k}{2}
\end{array}
\right] \right) +\left[ 
\begin{array}{l}
\,m \\ 
\frac{m}{2}
\end{array}
\right] ,$

\strut 

(3) if $m=k_{1}n_{j}+k_{2},$ for $k_{1}\in N$ and $1\leq k_{2}<n_{j},$ then

\strut 

$\ \ \ \ \ \ \ \ \ tr\left( T_{j}^{m}\right) =\left( 2^{k_{1}}-\left[ 
\begin{array}{l}
\,k_{1} \\ 
\frac{k_{1}}{2}
\end{array}
\right] \right) +\left[ 
\begin{array}{l}
\,m \\ 
\frac{m}{2}
\end{array}
\right] .$

$\square $
\end{theorem}

\strut \strut

We will finish this chapter with the following remark \strut ;

\strut

\begin{remark}
Let $G=\,<X:\,R>$ be a finitely presented group with its generator set $X$ $=
$ $\{x_{1},$ ..., $x_{N}\}$ and the relation $R$ $=$ $\{r_{1},$ ..., $r_{M}\}
$ and let's fix a generator $x_{j}\in X$ and the corresponding block
operator $T_{j}$ $=$ $x_{j}$ $+$ $x_{j}^{-1}$. Suppose there exist $r_{t}\in
R$ and $n_{j}\in \Bbb{N}$ such that $r_{t}=x_{j}^{n_{j}}.$ Then, by (2) and
(3) of the previous theorem,

\strut 

$\ \ \ \ tr\left( T_{j}^{m}\right) =\left\{ 
\begin{array}{ll}
\left[ 
\begin{array}{l}
\,m \\ 
\frac{m}{2}
\end{array}
\right]  & \text{if }m<n_{j} \\ 
&  \\ 
\left( 2^{k_{1}}-\left[ 
\begin{array}{l}
\,k_{1} \\ 
\frac{k_{1}}{2}
\end{array}
\right] \right) +\left[ 
\begin{array}{l}
\,m \\ 
\frac{m}{2}
\end{array}
\right]  & \text{if }m\geq n_{j},
\end{array}
\right. $

\strut 

where $m=k_{1}n_{j}+k_{2},$ for $k_{1}\in \Bbb{N}$ and $k_{2}\in \Bbb{N}\cup
\{0\}.$ $\ \square $
\end{remark}

\strut \strut

\strut

\strut

\subsection{Identically Distributedness}

\strut

\strut

By the previous section, we have that if $G$ $=$ $<X:\,R>$ is a finitely
presented group and if $x\in X,$ then

\strut

(i) \ if there exists $n\in \Bbb{N}$ such that $x^{n}=r,$ for $r\in R,$ then

\strut

\begin{center}
$tr\left( T^{m}\right) =\left\{ 
\begin{array}{ll}
\left[ 
\begin{array}{l}
\,m \\ 
\frac{m}{2}
\end{array}
\right] & \text{if }m<n_{j} \\ 
&  \\ 
\left( 2^{k_{1}}-\left[ 
\begin{array}{l}
\,k_{1} \\ 
\frac{k_{1}}{2}
\end{array}
\right] \right) +\left[ 
\begin{array}{l}
\,m \\ 
\frac{m}{2}
\end{array}
\right] & \text{if }m\geq n_{j},
\end{array}
\right. $
\end{center}

\strut

where $m=k_{1}n+k_{2},$ for $k_{1}\in \Bbb{N}$ and $k_{2}\in \Bbb{N}\cup
\{0\},$ and

\strut

(ii) if there is no $n\in \Bbb{N}$ such that $r=x^{n},$ for all $r\in R,$
then

\strut

\begin{center}
$tr(T^{m})=\left[ 
\begin{array}{l}
\,m \\ 
\frac{m}{2}
\end{array}
\right] ,$ \ for all $m\in \Bbb{N}.$
\end{center}

\strut

The above formuli directly proves the following theorem ;

\strut

\begin{theorem}
Let $G_{i}=\,<X_{i}:R_{i}>$ be finitely presented groups, for $i=1,2,$ and
assume that $x_{i}\in X_{i}$ are generators of $G_{i}$, for $i=1,2$. If
there are relators $r_{i}\in R_{i}$ and $n\in \Bbb{N}$ such that $r_{i}$ $=$ 
$x_{i}^{n},$ for all $i$ $=$ $1,$ $2,$ then the block operators $\left(
x_{i}+x_{i}^{-1}\right) $ of the group von Neumann algebras $L(G_{i}),$ $i$ $%
=$ $1,$ $2,$ are identically distributed.
\end{theorem}

\strut

\begin{proof}
By (i) and (ii) in the previous paragraph, case by case, we can get that

\strut

$\ \ \ \ \ \ \ \ \ tr(T_{1}^{m})=tr(T_{2}^{m}),$ \ for all $\ m\in \Bbb{N}.$

\strut
\end{proof}

\strut

In the above theorem, $G_{1}$ and $G_{2}$ are not necessarily distinct.
i.e., suppose that we have a finitely presented group $G=\,<X:\,R>$ and its
generators $x_{1},$ $x_{2}$ in $X$ satisfying that $x_{1}^{n}$ and $%
x_{2}^{n} $ are relators in $R.$ Then the block operators $\left(
x_{1}+x_{1}^{-1}\right) $ and $\left( x_{2}+x_{2}^{-1}\right) $ are
identically distributed in the group von Neumann algebra, $L(G).$

\strut

The following theorem is proved, similarly ;

\strut

\begin{theorem}
Let $G_{i}=\,<X_{i}:R_{i}>$ be finitely presented groups, for $i=1,2,$ and
assume that $x_{i}\in X_{i}$ are generators of $G_{i}$, for $i=1,2$. Suppose
that there is no numbers $n_{1}$ , $n_{2}$ $\in $ $\Bbb{N}$ such that $%
x_{i}^{n_{i}}$ $\in $ $R_{i},$ for $i$ $=$ $1,$ $2,$ then the block
operators $x_{i}+x_{i}^{-1}$ in $L(G_{i})$ are identically distributed. $%
\square $
\end{theorem}

\strut \strut \strut

\strut \strut \strut

\strut \strut \strut \strut \strut \strut

\strut \strut \strut

\strut \textbf{Reference}

\strut

\strut

\begin{quote}
{\small [1] \ \ A. Nica, R-transform in Free Probability, IHP course note,
available at www.math.uwaterloo.ca/\symbol{126}anica.}

{\small [2]\strut \ \ \ A. Nica and R. Speicher, R-diagonal Pair-A Common
Approach to Haar Unitaries and Circular Elements, (1995), www
.mast.queensu.ca/\symbol{126}speicher.\strut }

{\small [3] \ }$\ ${\small B. Solel, You can see the arrows in a Quiver
Operator Algebras, (2000), preprint.}

{\small [4] \ \ D. Shlyakhtenko, Some Applications of Freeness with
Amalgamation, J. Reine Angew. Math, 500 (1998), 191-212.\strut }

{\small [5] \ \ D.Voiculescu, K. Dykemma and A. Nica, Free Random Variables,
CRM Monograph Series Vol 1 (1992).\strut }

{\small [6] \ \ D. Voiculescu, Operations on Certain Non-commuting
Operator-Valued Random Variables, Ast\'{e}risque, 232 (1995), 243-275.\strut 
}

{\small [7] \ \ D. Shlyakhtenko, A-Valued Semicircular Systems, J. of Funct
Anal, 166 (1999), 1-47.\strut }

{\small [8] \ \ D.W. Kribs and M.T. Jury, Ideal Structure in Free
Semigroupoid Algebras from Directed Graphs, preprint.}

{\small [9] \ \ D.W. Kribs and S.C. Power, Free Semigroupoid Algebras,
preprint.}

{\small [10] I. Cho, Graph }$W^{*}$-{\small Probability Theory, (2004),
Preprint.}

{\small [11] I. Cho, Random Variables in Graph }$W^{*}$-{\small Probability
Spaces, (2004), Preprint. }

{\small [12] I. Cho, Amalgamated Semicircular Systems in Graph }$W^{*}$%
{\small -Probability Spaces, (2004), Preprint. }

{\small [13] I. Cho, Weighted Graph }$W^{*}${\small -Probability Spaces,
(2004), Preprint. }

{\small [14] I. Cho, The Moment Series of the Generating Operator of }$%
L(F_{2})*_{L(F_{1})}L(F_{2})${\small , (2003), Preprint. }

{\small [15] I. Grossman and W. Magnus, Groups and Their Graphs, MAA (1964),
ISBN-0-88385-600-X.}

{\small [16] K. J. Horadam, The Word Problem and Related Results for Graph
Product Groups, Proc. AMS, vol. 82, No 2, (1981) 157-164. }

{\small [17] J. Lauri and R. Scapellato, Topics in Graph Automorphisms and
Reconstruction, London Math Soc. Student No. 54, (1995), Cambridge Univ.
Press. }

{\small [18] P. Biane and R. Speicher, Stochastic Calculus with respect to
Free Brownian Motion and Analysis on Wigner Space, Prob. Theory Relat.
Fields 112, (1998) 378-409.}

{\small [19] R. Speicher, Combinatorial Theory of the Free Product with
Amalgamation and Operator-Valued Free Probability Theory, AMS Mem, Vol 132 ,
Num 627 , (1998).}

{\small [20] R. Speicher, Combinatorics of Free Probability Theory IHP
course note, available at www.mast.queensu.ca/\symbol{126}speicher.\strut }

{\small [21] S. Hermiller and J. Meier, Algorithms and Geometry for Graph
Products of Groups, J of Algebra 171, (1995) 230-257.}
\end{quote}

\end{document}